\documentclass{article} 
\usepackage{iclr2022_conference,times}

\usepackage{amsmath,amsfonts,bm}

\def\eqref#1{equation~\ref{#1}}

\def\1{\bm{1}}

\def\va{{\bm{a}}}
\def\vb{{\bm{b}}}

\def\ve{{\bm{e}}}

\def\vh{{\bm{h}}}

\def\vk{{\bm{k}}}

\def\vo{{\bm{o}}}
\def\vp{{\bm{p}}}
\def\vq{{\bm{q}}}

\def\vs{{\bm{s}}}

\def\vu{{\bm{u}}}
\def\vv{{\bm{v}}}

\def\vx{{\bm{x}}}

\DeclareMathAlphabet{\mathsfit}{\encodingdefault}{\sfdefault}{m}{sl}
\SetMathAlphabet{\mathsfit}{bold}{\encodingdefault}{\sfdefault}{bx}{n}

\usepackage{hyperref}
\usepackage{url}

\usepackage[utf8]{inputenc}
\usepackage{booktabs}
\usepackage[frozencache,cachedir=.]{minted} 
\usepackage{graphicx}
\usepackage{caption}
\usepackage{subcaption}
\usepackage{multirow}
\usepackage{xparse}
\title{RLOR: A Flexible Framework of Deep Reinforcement Learning for Operation Research}

\author{Ching Pui WAN\\
Lenovo Machine Intelligence Center\\Hong Kong\\
\texttt{cpwandeep@gmail.com}
\And
Tung LI\\
Lenovo Machine Intelligence Center\\Hong Kong\\
\texttt{tonyli537@gmail.com}
\And
Jason Min WANG\\
Lenovo Machine Intelligence Center\\Hong Kong\\
\texttt{jasonwangm@connect.ust.hk}
}

\NewDocumentCommand{\Colorbox}{O{\dimexpr\linewidth-2\fboxsep} m m}{%
  \colorbox{#2}{\strut \makebox[#1][l]{#3}}}

\setminted[python]{mathescape, escapeinside=@@,
               numbersep=5pt,
               gobble=2,
               frame=lines,
               framesep=2mm}
\definecolor{diffplus}{RGB}{230, 255, 236}
\definecolor{diffminus}{RGB}{255, 235, 233}

\iclrfinalcopy
\begin{document}

\maketitle
\begin{abstract}
Reinforcement learning has been applied in operation research and has shown promise in solving large combinatorial optimization problems. However, existing works focus on developing neural network architectures for certain problems. These works lack the flexibility to incorporate recent advances in reinforcement learning, as well as the flexibility of customizing model architectures for operation research problems. In this work, we analyze the end-to-end autoregressive models for vehicle routing problems and show that these models can benefit from the recent advances in reinforcement learning with a careful re-implementation of the model architecture. In particular, we re-implemented the Attention Model and trained it with Proximal Policy Optimization in CleanRL, showing at least 8 times speed up in training time. We hereby introduce RLOR, a flexible framework for Deep Reinforcement Learning for Operation Research. We believe that a flexible framework is key to developing deep reinforcement learning models for operation research problems. The code of our work is publicly available at \url{https://github.com/cpwan/RLOR}.
\end{abstract}

\section{Introduction}
Pointer Network \citep{Vinyals2015PointerNetworks} is a milestone work of applying neural networks in combinatorial optimization problems. It enabled dynamic input size and permutation invariance of input in the neural networks. In other words, we can feed a \textit{set} to a neural network. PN+RL \citep{Bello2019NeuralLearning} is another milestone work. It enabled training neural networks with reinforcement learning (RL) with REINFORCE algorithm, instead of requiring expensive ground truths from solvers for supervised learning. Since then, the REINFORCE algorithm (but rarely other RL algorithms) has been used in subsequent works for vehicle routing problems, including Order-invariant PN+ RL \citep{Nazari2018ReinforcementProblem}, Attention Model \citep{Kool2019AttentionProblems}, and POMO \citep{Kwon2020POMO:Learning}. Stemming from the same milestone works, DQN was preferred for graph problems, such as in S2V-DQN \citep{Dai2017LearningGraphs}, ECO-DQN \citep{Barrett2020ExploratoryLearning}, and GCOMB \citep{Manchanda2020GCOMB:Graphs}. These streams of work focused on improving model architectures instead of exploring more advanced RL algorithms.

On the other hand, the RL community has been growing. Several RL platforms have been proposed for academic research and industrial applications, including StableBaselines3 \citep{stable-baselines3}, RLlib\citep{Liang2017RLlib:Learning}, DI-Engine\citep{ding}, Tianshou\citep{tianshou}, and CleanRL\citep{Huang2022CleanRL:Algorithms}. These RL platforms have their design philosophies and different levels of abstraction to accommodate numerous RL algorithms. Nevertheless, they share several similarities. In contrast to the complicated model architectures used for vehicle routing problems, an MLP model architecture and a 1-d vector input are assumed in most of the RL platforms. Moreover, REINFORCE (or policy gradient), the most popular algorithm for end-to-end vehicle routing problems, was not implemented in most of these RL platforms. It was claimed by their developers that the 31 years old REINFORCE algorithm \citep{Williams1992SimpleLearning} usually does not perform well compared to the recent RL algorithms. The disentanglement has been demonstrated between the advances in model architecture and RL algorithms.

A natural question arises: can we leverage the advanced algorithms in the RL platforms for vehicle routing problems? In this work, we will demonstrate a \textbf{Yes} to this question. However, it is a challenging task, due to two major reasons: compatibility and efficiency. In the RL platforms, they accept agents built with an MLP model, a CNN model, or even an RNN model, but none of them has considered the attention model in their design. (DI-Engine \citep{ding} implements DecisionTransformer \citep{Chen2021DecisionModeling} but it is considered as an algorithm instead of an agent building block) As a result, it would require a huge effort to modify the RL platform to allow nested observations, dynamic input size, and management between hidden states. Eventually, the model architecture can be trained in the RL platform after fixing the compatibility issue. However, it still suffers from the efficiency issue. For instance, the RL platforms have their own environment APIs and there are different levels of overhead in data transformation and device communications (between CPU/GPU). We performed experiments on several RL platforms and find that CleanRL \citep{Huang2022CleanRL:Algorithms} has the lowest overhead and thus the highest efficiency.

We introduce the RLOR framework, which consists of four parts: \textit{model}, \textit{algorithm}, \textit{environment}, and \textit{search}. The \textit{model} describes the neural network model architecture. Our default model is developed and refactored from the Attention Model \citep{Kool2019AttentionProblems}. The \textit{algorithm} describes the RL algorithm. We adapted the Proximal Policy Optimization (PPO) algorithm from CleanRL to our problems. The \textit{environment} describes the RL environments. We followed the protocol of the OpenAI Gym \citep{1606.01540} when defining RL environments for the operation research problems. The \textit{search} defines the decoding strategies of the neural network, such as greedy or sampling.

Our major contributions are in three-folds:
\begin{enumerate}
    \item As far as we know, it is the first work to incorporate end-to-end vehicle routing model in modern RL platforms\\
    \item It speeds up the training of Attention Model by 8 times (e.g., 25 hours $\to$ 3 hours)\\
    \item It introduces a flexible framework for developing \textit{model}, \textit{algorithm}, \textit{environment}, and \textit{search} for operation research
\end{enumerate}

In Section \ref{sec:relatedworks}, we will discuss the related works. In Section \ref{sec:model}, \ref{sec:algo},\ref{sec:env}, \ref{sec:search}, we will discuss how we solve the compatibility issues and integrate \textit{model}, \textit{algorithm}, \textit{environment}, and \textit{search} for OR models in RL platforms. In Section \ref{sec:efficiency}, we will discuss how we solve the efficiency issues of the OR model in RL platforms. In Section \ref{sec:results}, we will discuss the results of our experiments. We will refer ``model" as the neural network and ``OR models" as neural networks applied to operation research problems. We will refer ``trajectory" as the sequence of actions and states until a given number of steps while a ``rollout" is a complete trajectory.

\section{Related Work}\label{sec:relatedworks}
\paragraph{Supervised learning methods.}
Earlier deep learning approaches to solving combinatorial optimization problems mostly featured supervised learning. Pointer Network \citep{Vinyals2015PointerNetworks} formulated the combinatorial optimization problem as a sequence-to-sequence problem and predicted the target sequence obtained from a numerical solver. Later works \citep{joshi2019efficient, Li2018CombinatorialSearch} employed Graph Neural Network to predict the adjacency matrix of the optimal solution and performed searches over the adjacency matrix, showing better performance. However, their supervised learning natures implied expensive training set collection. The resulting model also had limited generalization power. Generalized GNN \citep{fu2020GeneralizeGNN} and DPDP \citep{kool2021DynamicProgramming} tried to address the generalization issues with GNN model by introducing advanced search techniques (Monte Carlo Tree Search, dynamic programming). However, the analysis from \citet{bother2022s} suggested that the performance of the supervised learning model may (in general) be contributed by the search routine instead of the model architecture. Therefore, we turn our attention to the deep RL approach.

\paragraph{Deep RL methods.}
There are two streams of approaches for applying RL in solving combinatorial optimization problems: the construction method, and the improvement method. Construction methods, such as PN+RL \citep{Bello2019NeuralLearning}, Attention Model \citep{Kool2019AttentionProblems}, POMO \citep{Kwon2020POMO:Learning}, constructed the solution step by step. On the other hand, the improvement methods such as Learning Improvement Heuristics \citep{wu2019LearningImprovementHeuristics},  \citep{Lu2020AProblems}, NeuRewritter \citep{Chen2019LearningOptimization}, DACT \citep{Ma2021LearningTransformer}, started with a complete solution and predicted the location for local rewriting to search for a better solution. The improvement methods were much slower than the construction methods but could obtain a better solution. Instead of predicting the local rewriting with a neural network, eMagic \citep{ouyang2021EquivarianceGNN} performed local search (a heuristic) to find better rollouts during the RL training of their construction method. There were also efforts to combine tree search techniques from constraint programming in the work of \citet{Cappart2021CombiningOptimization}. Readers interested in the recent advances in neural combinatorial optimization can refer to the survey blog of \citet{chaitanyak2022recentadvancesin}. 

\paragraph{RL libraries.}
RL platforms provide training pipelines for a variety of RL algorithms. Popular RL platforms include StableBaselines3 \citep{stable-baselines3}, RLlib \citep{Liang2017RLlib:Learning}, DI-Engine \citep{ding}, Tianshou \citep{tianshou}, and CleanRL \citep{Huang2022CleanRL:Algorithms}. In our work, we employ the PPO algorithm provided by these RL platforms. Readers interested in the implementation of PPO can refer to the blog of \citet{shengyi2022the37implementation}. On the other hand, we reformulated OR problems into RL environments so that they can communicate with RL algorithms through the environment API. OpenAI Gym \citep{1606.01540} provides standardized API defining RL environments and is supported in most of the RL platforms. OR-Gym \citep{HubbsOR-Gym} implements RL environments for a couple of small-scale operation research problems. Graphenv \citep{biagioni2022graphenv} implements the RL environment for travelling salesman problem (TSP) while conforming to RLlib's environment API.

\section{Model}\label{sec:model}
\begin{figure}
    \centering
    \includegraphics[width=\textwidth]{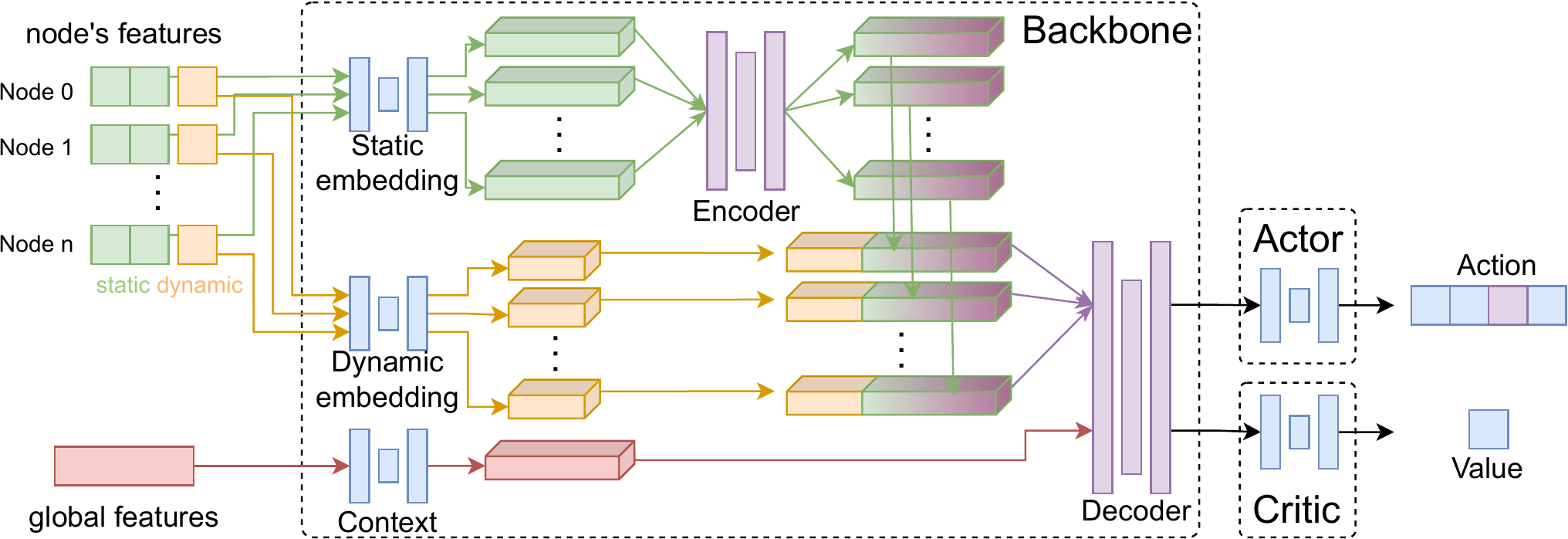}
    \caption{The model architecture (Backbone, Actor, and Critic) of the refactored Attention Model. The neural networks in blue represent MLPs while those in purple represent multi-head attention layers.}
    \label{fig:refactor-am}
\end{figure}
We performed detailed code reviews for the Pointer Network \citep{Vinyals2015PointerNetworks} and the Attention Model \citep{Kool2019AttentionProblems}. The pseudocodes for these model architectures can be found in Appendix \ref{appendix-models}. We refactored both models and upgraded them to the latest PyTorch 1.13. In this section, we will focus on the Attention Model, which we found easier to integrate into RL platforms. 

The overview of the refactored Attention Model can be found in Figure \ref{fig:refactor-am}. We will use the capacitated vehicle routing problem (CVRP) as an example to illustrate the workflow of the refactored Attention Model. The node-wise static features (coordinates of customers and depots) and the node-wise dynamic features (demands of customers) are projected to the static embedding and dynamic embedding independently for each node. The static embeddings of different nodes are integrated in the encoder. The dynamic embedding and the encoder output are concatenated and serve the key and value of the decoder. The global features (the current load of the vehicle) are projected to a latent representation of the global context and serve the query of the decoder. The decoder predicts a probability distribution of the next action (which node to visit next) along with the latent representation of the current environment state. The original implementation of Attention Model \citet{Kool2019AttentionProblems} concluded the model at the decoder. In our refactored version, we add additional modules --- Actor and Critic to better conform to the design of RL algorithms.

\textbf{Backbone.} We regard the aforementioned steps from input embedding to the decoder as the backbone. It is responsible for most of the computation. The backbone provides the probability distribution for the actor model and latent representation of the environment state to the critic model.

\textbf{Actor model.} The actor model selects an action from the probability distribution and sends it to the RL environment. The RL environment then returns the next observation. The state (node's features and global features) is updated accordingly and is used for the next iteration.

\textbf{Critic model.} The critic model predicts the \textit{value} of the current state (i.e., the average reward we can claim starting from the current state). To predict the value effectively, we need to feed the critic model with sufficient information about the state. In particular, we feed the critic model with the glimpse of the last attention layers of the decoder.

Our reformulation of the Attention Model provides a more compatible interaction among the neural network model, the RL algorithm, and the RL environment.

\subsection{Adapting OR models to RL platforms}\label{sec:adapt}
\begin{figure}
     \centering
     \begin{subfigure}[b]{0.48\textwidth}
         \centering
         \includegraphics[width=\textwidth]{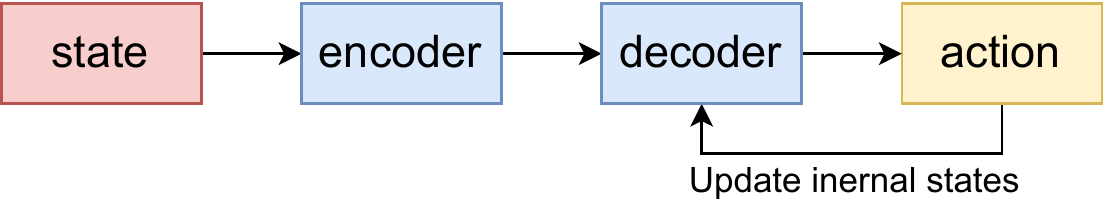}
         \caption{OR model pipeline}
         \label{fig:or model flow}
     \end{subfigure}
     \hfill
     \begin{subfigure}[b]{0.48\textwidth}
         \centering
         \includegraphics[width=\textwidth]{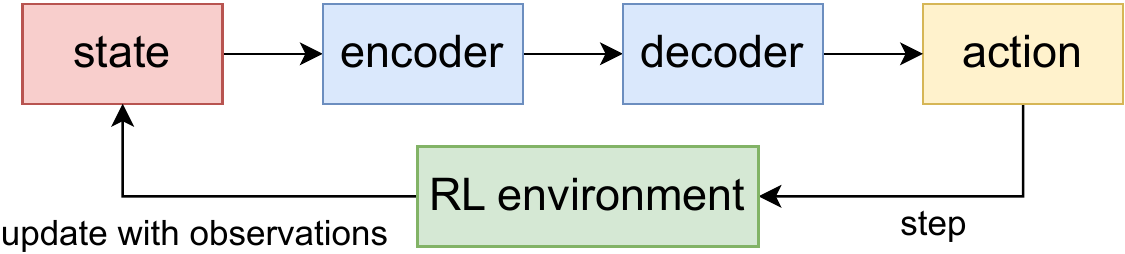}
         \caption{RL platform pipeline}
         \label{fig:rl platform flow}
     \end{subfigure}
     \caption{Comparison of action sampling pipelines between OR model and RL platform. The OR model samples action from the decoder and updates the internal states directly until all steps in a rollout are completed. In contrast, RL platforms sample action step-by-step during the interaction with the RL environment and require a complete forward pass of the neural network.}
     \label{fig:compare-pipelines}
\end{figure}

Since the RL platforms employ step-wise sampling, instead of rollout-wise sampling in OR models, we need to modify the OR models. In Figure \ref{fig:compare-pipelines}, we illustrate the difference in the pipelines. In RL platforms, a forward pass of the entire neural network model is performed and only the action can be used to update the next state. This is in contrast to OR models, where the internal state of the neural network can also affect the next prediction. 

\textbf{Pointer Network.}
The step-wise sampling of environment interaction makes it difficult to implement Pointer Network in RL platforms. In Pointer Network, the hidden state of the LSTM layers in the decoder is used to pass information across steps in a rollout while the output of the encoder is reused for each of the decoding steps. Some RL platforms (RLlib, Tianshou) support RNN models by passing the detached hidden state as additional information along with the observations, allowing the decoder to be trained. However, if we also treat the encoder output as additional information, the gradient would not be back-propagated to the encoder. Alternatively, some RL platforms (Tianshou) allow passing a non-detached hidden state, but they do not guarantee that hidden states can be consumed in the same iteration. As a result, the non-detached encoder output may be consumed in later iterations when the model has already been updated, resulting in a runtime error in PyTorch.

\textbf{Attention Model.}
The same difficulty exists in implementing the Attention Model in the RL platforms. However, there is a workaround --- with a sacrifice of efficiency. In our refactored Attention Model, the decoder does not depend on its previous hidden states. Instead, it can obtain all step-wise information from the state extracted from the RL environment observation. Hence, for each RL environment step, we can perform a feed-forward for the entire refactored Attention Model as shown in Figure \ref{fig:rl platform flow}. Each component of the Attention Model can be trained with this paradigm. However, the same encoder output is recomputed every time. In most of the RL platforms, it is unavoidable to make this trade-off in order to enable the training of the Attention Model.

\section{Algorithm}\label{sec:algo}
\subsection{RL algorithms in OR models}
In preceding works for vehicle routing problems, REINFORCE algorithm was used. In RL platforms, the REINFORCE algorithm is more often referred as Policy Gradient (PG). The goal of policy gradient with respect to an environment instance $s$ is to maximize the expected return
\begin{equation}\label{eq-reinforce}
    J(\theta | s) =\mathbb{E}_{\pi \sim p_{\theta}(.|s)}[L(\pi|s)]
\end{equation}
where $\theta$ is the neural network parameter, $\pi$ is the policy (in OR model context, the complete sequence of actions) drawn from the actor $p_{\theta}(.|s)$, and $L(\pi|s)$ is the return of executing policy $\pi$ on the environment instance $s$. 

The REINFORCE algorithm derived the gradient of the expected return as the following:
\begin{equation}\label{eq-reinforcegrad}
    \nabla_{\theta}J(\theta | s) = \mathbb{E}_{\pi \sim p_{\theta}(.|s)} [(L(\pi |s)-b(s)) \nabla_{\theta}\log(p_{\theta}(\pi|s))].
\end{equation}
where $b(s)$ is the baseline. Intuitively, the return $L(\pi|s)$ is compared against the baseline $b(s)$. The better (relatively) the return $L(\pi|s)$, the higher the log-likelihood the policy $\log(p_{\theta}(\pi|s))$ contributes to the gradient. Theoretically, a good choice of baseline $b(s)$ reduces the variance of the gradient and potentially leads to faster convergence. The preceding works for vehicle routing problems designed different baselines for the policy gradient. Interested readers may refer to Appendix \ref{appendix-algo} for the modifications they made.
\subsection{RL algorithms in RL platforms}
The formulation of policy gradient and its variants in RL platforms follows a similar derivation from the REINFORCE algorithm. However, the policy $\pi$ in Equation \ref{eq-reinforce} is explicitly broken down into a sequence of actions $\pi = \pi_0, \pi_1, ... , \pi_T$ and the expectation in Equation \ref{eq-reinforcegrad} is then with respect to each of these actions. The return $L(\pi |s)$ is broken down to the sum of rewards $r(\pi_t |s_t)$. The $L(\pi |s)-b(s)$ term is replaced by a step-wise estimate $L(\pi_t |s_t)-b(s_t)$, in which $L(\pi_t |s_t)$ is the sum of future rewards and $b(s_t)$ can be predicted by the critic. As an illustration of the algorithmic difference, we may compare Equation \ref{eq-reinforcegrad} with the gradient computed from the step-wise estimate in RL platforms:

\begin{equation}\label{eq-A2Cgrad}
    \nabla_{\theta}J(\theta | s) = \mathbb{E}_{t\sim \{0...T\}, \pi_t \sim p_{\theta}(.|s_t)} [(L(\pi_t |s_t)-b(s_t))\nabla_{\theta}\log(p_{\theta}(\pi_t|s_t))].
\end{equation}

Variants of policy gradient use different step-wise estimate. For example, A2C uses the advantage $A(\pi_t|s_t)$ computed from the rewards and critic predictions; PPO further applies clipping on the advantage with a regulation on the deviation of the policy. The comparison of RL algorithms in OR models and RL platforms is summarized in Table \ref{tab:comparison-in-RL-algo}. We accommodated the pipeline of the OR model to that of the RL platforms in Section \ref{sec:adapt}.

\begin{table}[t]
\begin{center}\resizebox{\textwidth}{!}{
\begin{tabular}{lll}
\toprule
                   & OR models                                            & RL platforms                                   \\ \midrule 
Algoritms & REINFORCE & A2C, PPO, DQN, ...\\
Unit of sample     & a full rollout $\{s_0,\pi_0,s_1,...,\pi_T,s_{T+1}\}$ & one step $s_t, \pi_t,s_{t+1}$ \\
Actor predicts...  & $n$ actions in one go                                & $1$ action each time                           \\ 
Critic predicts... & not used, baseline is derived from previous rollouts & the value of the current state                 \\ \bottomrule
\end{tabular}
}
\end{center}
\caption{Comparison of the RL algorithm used in OR models and RL platforms.}
\label{tab:comparison-in-RL-algo}
\end{table}

\section{Environment}\label{sec:env}

The primary settings and constraints of the vehicle routing problem environments are referenced from Attention Model paper \citep{Kool2019AttentionProblems}, and adapted into OpenAI Gym \citep{1606.01540}. The OpenAI Gym environment interface is supported in most RL platforms. The environment for each problem needs to be carefully designed such that the necessary information would be included in the observation returned by the OpenAI Gym environment, and a wrapper is required to transform the observation into the schema that the RL platforms intend to receive.

\subsection{Interface}

The following steps of the OR model involve interactions with the environment: 
\begin{itemize}
    \item Initialize the environment for the problem instance $s$. 
    \item Update the environment with the selected action $\pmb{\pi}_t$.
    \item Get whether the policy $\{\pmb{\pi}_t\}$ has finished.
    \item $\mathrm{Context}$ needs to get the global information, such as remaining vehicle capacity.
    \item $\mathrm{DynamicEmbedding}$ needs to get the dynamic node features.
    \item $\mathrm{mask_t}$ needs to get the feasible actions.
\end{itemize}

An RL environment for vehicle routing problems should provide the following interfaces:
\begin{minted}{python}
class RequiredEnv():
    def __init__(self, s: problem_instance): ...
    def step(self, action): ...
    def get_global_context(self): ...
    def get_dynamic_node_features(self): ...
    def get_mask(self): ...
    def all_finished(self): ...
\end{minted}

Compare to the OpenAI Gym environment interface:
\begin{minted}{python}
class GymEnv():
    def __init__(self, arg1, arg2, ...):
    # What action can be taken?
    self.action_space = ...
    # What information can we get from the problem?
    self.observation_space = ...
    def step(self, action):
        ...
        return observation, reward, done, info
    def reset(self):
        ...
        return observation  
\end{minted}

To translate the \mintinline{python}{GymEnv} to our \mintinline{python}{RequiredEnv}, we need to store the results of the \mintinline{python}{GymEnv.step(action)}. When the \mintinline{python}{RequiredEnv.get_xxxx} function is called, we extract the corresponding information from the stored \mintinline{python}{observation}. When the  \mintinline{python}{RequiredEnv.all_finished} function is called, we extract the corresponding information from the stored \mintinline{python}{done}.

\subsection{Adaptation to RL libraries}
Here, we give an example to adapt the RL libraries interface to the Attention Model. We use CleanRL as an example for demonstration purpose.
In CleanRL, the observations from the environment is referred as \mintinline{python}{next_obs} and is passed to the model in each step of an episode. The \mintinline{python}{next_obs} is a \mintinline{python}{dict} of tensors. In particular, for TSP, the \mintinline{python}{next_obs} has the following contents:

\begin{center}
\begin{tabular}{cccl}\toprule
\mintinline{python}{next_obs.keys()} & \mintinline{python}{shape} & \mintinline{python}{type} & descriptions \\
\midrule
\mintinline{python}{'observations'}& \mintinline{python}{[B,n,2]}      &\mintinline{python}{float} &   coordinates of nodes          \\
\mintinline{python}{'action_mask'}& \mintinline{python}{[B,n] }     & \mintinline{python}{float} &  0: forbidden; 1: permitted          \\
\mintinline{python}{'first_node_idx'}&  \mintinline{python}{[B] }    & \mintinline{python}{int} & the index of the first node visited       \\
\mintinline{python}{'last_node_idx'}& \mintinline{python}{[B]}      & \mintinline{python}{int} & the index of the previous node visited       \\
\mintinline{python}{'is_initial_action'}&   \mintinline{python}{[B]}    &  \mintinline{python}{bool} &  whether the current step is the first step     \\
\bottomrule
\end{tabular}
\end{center}

On the other hand, the Attention Model for TSP requires a couple of data members and functions from a \mintinline{python}{state} object:

\begin{center}
\begin{tabular}{cccl}\toprule
Required data/ functions & \mintinline{python}{shape} & \mintinline{python}{type} & descriptions \\
\midrule
\mintinline{python}{state.get_mask()}& \mintinline{python}{[B,1,n] }     & \mintinline{python}{bool} &  \mintinline{python}{True}: forbidden, \mintinline{python}{False}: permitted          \\
\mintinline{python}{state.first_a}&  \mintinline{python}{[B,1] }    & \mintinline{python}{int} & the index of the first node visited       \\
\mintinline{python}{state.get_current_node()}& \mintinline{python}{[B,1]}      & \mintinline{python}{int} & the index of the previous node visited       \\
\mintinline{python}{state.is_initial_action}&   \mintinline{python}{[B]}    &  \mintinline{python}{bool} &  whether the current step is the first step     \\
\bottomrule
\end{tabular}
\end{center}

Hence, we define a wrapper to translate from the nested observations dict \mintinline{python}{(next_obs : dict)} collected by CleanRL to the object \mintinline{python}{(state : object)} that our model required.
\begin{minted}{python}
# wrapper for the problem
class stateWrapper:
    """
    From dict of numpy arrays 
        to an object that supplies PyTorch tensors.
    """

    def __init__(self, next_obs, device, problem):
        self.device = device
        self.states = {k: torch.tensor(v, device=self.device) 
                        for k, v in next_obs.items()}
        if problem == 'tsp':
            self.is_initial_action = \
                self.next_obs["is_initial_action"].to(torch.bool)
            self.first_a = self.next_obs["first_node_idx"]
        elif problem == 'cvrp':
            input = {'loc': self.next_obs['observations'], 
                    'depot': self.next_obs['depot'].squeeze(-1),
                    'demand': self.next_obs['demand']}
            self.states['observations'] = input
            self.VEHICLE_CAPACITY = 0
            self.used_capacity = -self.next_obs["current_load"]

    def get_current_node(self):
        return self.states["last_node_idx"]

    def get_mask(self):
        return (1 - self.states["action_mask"]).to(torch.bool)

# Inside the Attention Model
state = stateWrapper(next_obs)
\end{minted}

\section{Search}\label{sec:search}
During inference, we are searching for a solution with the probability distribution predicted by the Attention Model. In RL platforms, the action of an agent is mostly selected in two ways: greedy (choosing the highest probability action), or sampling (selecting an action according to the probability distribution). In contrast, some works on OR models borrow beam search from NLP \citep{Nazari2018ReinforcementProblem} or fine-tune the model on test instance with Efficient Active Search \citep{Hottung2022EfficientProblems}. In the work by \citet{Cappart2021CombiningOptimization}, they treated the probability distribution as a heuristic and performed a tree search accordingly. In our implementation, we borrow the idea of POMO \citep{Kwon2020POMO:Learning} for searching. In POMO, it performed greedy rollouts with different starting nodes during inference. We denote this decoding strategy as Multi-Greedy.

\section{Efficient Training on RL platform}\label{sec:efficiency}
We implemented our refactored Attention Model in 4 major RL platforms: DI-engine, RLlib, Tianshou, and CleanRL. We did not implement it on StableBaselines3 due to engineering difficulties. We found that CleanRL has the highest efficiency in environment step collection and the best convergence. The results will be discussed in Section \ref{sec:results}. We found that even the best RL platform could not beat the performance of the original implementation. We performed a detailed analysis and developed several ways to improve the training efficiency of our OR model on CleanRL.
\subsection{Reducing communication overhead}
We found that the communication overhead between GPU/ CPU devices is one of the major overheads in the RL platforms. For example, we profiled DI-engine and found that 30\% of the computation was spent on data manipulation for observations and moving the data across devices. In CleanRL, we implement the RL environment with Numpy and convert the observations (in NumPy array) to PyTorch Tensor on GPU in the state wrapper. That is, data collected from the environment are moved to GPU only when it is needed for the forward pass.  The ad-hoc conversion of data across devices reduces the communication overhead and makes CleanRL the most efficient RL platform for our problems.
\subsection{Cached Encoder}\label{sec:cached-encoder}
In Section \ref{sec:adapt}, we mentioned that we have to sacrifice the efficiency of the OR model to accommodate the workflow of the RL platforms. To retain the efficiency, we modify the CleanRL to allow caching of the encoder output. We enable cached encoder during both value bootstrapping and model training. The encoder output is computed once and is passed to the decoder for sampling the action and values until the environment is done. We enforce that the encoder output has to be consumed in the same iteration. The decoder is fed with the non-detached version of the encoder output. The back-propagation is performed only when the forward pass for every step is done. In this way, the encoder can be updated with the gradients computed from every environment step. We profiled the implementation and found that it could save $\sim$ 80\% computation during forward pass by caching the encoder output.

\subsection{Larger Batch}
With the PPO algorithm, the training is more stable and we can use a larger learning rate with a larger batch-size. In particular, we collect $1024\cdot 51$ steps in a batch and train our model with 8 minibatch (minibatch-size = 6528) for the 50-nodes TSP. A larger learning rate of $10^{-3}$ can be used without weight decay. It takes 6.9 GiB GPU memory. It may scale further but we want to keep it compact so that it can be fitted in an affordable GPU (for instance, those in the free-tier Colab notebook).

\subsection{Parallel decoding}
In Section \ref{sec:cached-encoder}, we share the encoder output across environment steps. In this section, we further speed up computation by sharing the encoder output across different trajectories of the same RL environment instance. The dynamic nodes features and global context constitute the query in the attention layers of the decoder. The Attention Model only passes one query (which corresponds to one trajectory) to the decoder per environment. We can pass multiple queries, each corresponding to different trajectories, to the decoder. By decoding multiple trajectories in parallel in the GPU, we obtain better efficiency than decoding them separately. The pseudo-code of the attention layer can be found in Appendix \ref{appendix-attention-model}. We replace the query $q$ with a higher dimension tensor $\vq$. Note that, during training, we sample the trajectories according to the probability distribution predicted by the decoder. This is in contrast to POMO \citep{Kwon2020POMO:Learning} in which different starting nodes are enforced.

\subsection{Environment Vectorization}

We implemented a bi-level vectorized environment. The first level is the OpenAI Gym vector environment consisting of multiple sub-environments. In each of these sub-environments, we implemented another layer of sub-environment --- vectorized environment written in Numpy's vectorized operations. For example, we can create a $M\times N$ bi-level vectorized environment, in which there are $M$ RL environment instances and $N$ trajectories for each of the instances. We profiled the implementation and found that the bi-level environment vectorization can bring a 10x speed up for environment steps collection for $M=1024,N=50$. More on the design of the bi-level environment vectorization can be found in Appendix \ref{appendix-vectorization}

\section{Experiment Results}\label{sec:results}

\begin{table}[t]
\begin{center}
  \begin{tabular}{l|ccrr|ccrr}
    \toprule
    \multirow{2}{*}{Method} &
      \multicolumn{4}{c}{TSP50} &
      \multicolumn{4}{c}{CVRP50}\\
    & Len. & Gap & Time & Steps & Len. & Gap & Time & Steps\\
    \midrule
    Optimal & 5.69 & 0.00\% & - & - & 10.38 & 0.00\% & - & -\\
    \midrule
    \citet{Kool2019AttentionProblems} (PG) & 5.80 & 1.93\% & 25h & 6.5G & 10.98 & 5.78\% & 33h & 7.8G\\
    \midrule
    DI-Engine (PPO) & 5.90 & 3.69\% & 3d & 480M & - & - & -  & - \\
    RLlib (PPO) & 6.15 & 8.08\% & 3d & 153M & - & - & - & - \\
    Tianshou (PPO) & 6.19 & 8.79\% & 3d & 611M & - & - & - & - \\
    CleanRL (PPO) & 5.83 & 2.46\% & 3d & 1.7G & - & - & - & - \\
    \midrule
    Ours (PPO) & 5.79 & 1.76\% & 24h & 30G & 10.91 & 5.11\% & 24h & 21G \\
    + Multi-Greedy & \textbf{5.71} & \textbf{0.35}\% & 24h & 30G & \textbf{10.82} & \textbf{4.23}\% & 24h & 21G \\
    \bottomrule
  \end{tabular}
\end{center}
\caption{Experiment results across implementations and RL platforms. The second row refers to the original implementation of the Attention Model. The middle rows refer to our implementation on different RL platforms. The last two rows refer to our implementation with a bag of tricks. The RL algorithm in each implementation is indicated in the bracket. Experiments are done on 50-nodes TSP and 50-nodes CVRP environment.}
\label{tab:results}
\end{table}

\begin{table}[t]
\begin{center}
  \begin{tabular}{l|r|r|c}
    \toprule
    Method & Time & Env Steps& Performance\\
    \midrule
    \citet{Kool2019AttentionProblems} (PG) & 50m & 65M & 6.00 \\
    \midrule
    DI-Engine (PPO) & 12h & 40M & 6.00\\
    RLlib (PPO) & - & - & 6.00\\
    Tianshou (PPO) & - & - & 6.00\\
    CleanRL (PPO) & 3h & 70M & 6.00\\
    \midrule
    Ours (PPO) & \textbf{14m }& 250M & 6.00\\
    \bottomrule
  \end{tabular}
\end{center}
\caption{Time and environment steps required to reach 6.0 in route length ($\sim$ 5\% gap) in the 50-nodes TSP. The first row refers to the original implementation of the Attention Model. The middle rows refer to our implementation on different RL platforms. The last row refers to our implementation with a bag of tricks. The RL algorithm in each implementation is indicated in the bracket. RLlib and Tianshou did not attain the performance in a reasonable time.}\label{tab:platforms}
\end{table}

\paragraph{Training.} We first describe our experiment settings. We used the refactored Attention Model as our agent for all OR problems. We trained our agent with the PPO algorithm in RL platforms. The training was done with the Adam optimizer with a learning rate of $\eta=10^{-3}$ on a single Titan RTX GPU. Each epoch consists of 1024 randomly generated RL environment instances. For comparison, we also trained the original implementation of the Attention Model from \citet{Kool2019AttentionProblems}. We followed its default setting and trained it for 100 epochs on the same GPU.

\paragraph{Inference.} For each problem, the performance is measured on the 10000 test instances as being done in \citet{Kool2019AttentionProblems}. The optimality gap is calculated with respect to the length of the optimal route of the vehicle routing problem instance. Greedy search is the default decoding strategy for the Attention Model. For Multi-Greedy search, we perform $n=50$ greedy searches with different start nodes and select the best trajectory.

\subsection{Attention Model On Different RL platforms}

We integrated our refactored Attention Model on several RL platforms and trained it with the PPO algorithm provided in these RL platforms. The results are summarized in Table \ref{tab:results}. We found that none of these RL platforms were able to attain the performance of the original implementation of the Attention Model. We speculate that the discrepancy between RL platforms and the original implementation originated from the efficiency of the model-environment interactions and data manipulations. In particular, CleanRL employed the most efficient vector API (with the fastest environment step collection among the RL platforms), and produced the best results among the RL platforms. We observed a positive correlation between the number of environment steps and the final performance across the RL platforms. It indicated that the efficiency in environment step collection could be vital for reaching a tighter optimality gap. Based on this assumption, we improved the training efficiency on CleanRL with a bag of tricks which has been elaborated on in Section \ref{sec:efficiency}. With the improved implementation, we were able to surpass the Attention Model \citep{Kool2019AttentionProblems} in both training time and the optimality gap. Our approach collected 5 times more environment steps and reached a better objective given at the same time. With Multi-Greedy search, we can obtain an even better performance. We observed a similar performance boost for the CVRP50 instances.

From another perspective, our approach can reach the same optimality gap with a shorter training time. Table \ref{tab:platforms} reported the time required to reach 5\% gap for the 50-nodes TSP. Among the others, our approach was the most efficient and it took only 14 minutes for the training. Our approach is appealing in that it can train a sufficiently good OR model in a short training time. It can accelerate the development of a new OR model. For instance, hyper-parameter tuning can be done in a shorter time.

\begin{table}[t]
\begin{center}
  \begin{tabular}{l|r|r|c}
    \toprule
    Method & Training Time & Env Steps & Performance\\
    \midrule
    \citet{Kool2019AttentionProblems} & 25h & 6.5G & 5.80\\
    \midrule
    Baseline & 11d & 6G & 5.79\\
    + Cached Encoder & 4d & 6G & 5.79\\
    + Larger Batch & 2.5d & 6G & 5.79\\
    + Parallel Decoding \& Env Vectorization & 24h & 32G & 5.79\\
    + Multi-Greedy & 3h & 5G & 5.79\\
    \bottomrule
  \end{tabular}
\end{center}
\caption{Ablation study on different components of our framework. The first row refers to the original implementation of Attention Model. The following rows refer to the designs in our framework. Experiments were done on the 50-nodes TSP.}
\label{tab:ablation}
\end{table}


\subsection{Ablation Study}

As an ablation study, we performed experiments with CleanRL with different parts of our design. The results are shown in Table \ref{tab:ablation}. We can observe significant improvements with the components. By caching the encoder outputs, we were able to greatly reduce unnecessary re-computation during training. We profiled the implementation and found that 80\% of the training time can be saved by caching the encoder output. By using a larger batch with a higher learning rate, we could scale the training while fitting the training inside a single GPU. By using parallel decoding and the bi-level environment vectorization, we were able to substantially increase the environment steps trained per second, boosting the collection efficiency of the environment observations. With these designs, our implementation has already surpassed the original implementation of the Attention Model. Lastly, by employing Multi-Greedy search during inference, we were able to \textbf{achieve the same performance in just 3 hours}.

\begin{table}[t]
\begin{center}
  \begin{tabular}{l|cccc|cccc}
    \toprule
    \multirow{2}{*}{Method} &
      \multicolumn{4}{c}{TSP50} &
      \multicolumn{4}{c}{CVRP50}\\
    & Len. & Gap & Time & Steps & Len. & Gap & Time & Steps\\
    \midrule
    Optimal & 5.69 & 0.00\% & - & - & 10.38 & 0.00\% & - & -\\
    \midrule
    POMO (PG) & 5.73 & 0.70\% & 24h & 61G & 10.84 & 4.43\% & 24h & 55G \\
    + Multi-Greedy & \textbf{5.71} & \textbf{0.32}\% & 24h & 61G & \textbf{10.60} & \textbf{2.21\% }& 24h & 55G \\
    \midrule
    Ours (PPO) & 5.79 & 1.76\% & 24h & 30G & 10.91 & 5.11\% & 24h & 21G \\
    + Multi-Greedy & 5.71 & 0.35\% & 24h & 30G & 10.82 & 4.23\% & 24h & 21G \\
    \bottomrule
  \end{tabular}
\end{center}
\caption{Experiment results compared to POMO. The RL algorithm in each implementation is indicated in the bracket. Experiments are done on 50-nodes TSP and 50-nodes CVRP.}
\label{tab:POMO}
\end{table}

\subsection{Potential extension with POMO}
As a comparison, we also trained the implementation from POMO \citep{Kwon2020POMO:Learning}, a state-of-the-art construction method. We used equivalent model settings (e.g., 3 encoder layers) and trained POMO for 24 hours. The result is reported in Table \ref{tab:POMO}. We note that our method did not outperform POMO. In our design, we attempt to modularize the \textit{model}, \textit{algorithm}, \textit{environment}, and \textit{search}, in order to create a flexible and comprehensive framework for developing deep RL model for the operation research problem. We made a trade-off between flexibility and the environment step collection efficiency. For instance, we could implement the entire vectorized environment in PyTorch and obtain higher efficiency in environment step collection. Yet, it will incur burdens on fine-grained control over the RL environment. Nevertheless, we are still close behind the state-of-the-art benchmark, with a clean one-page code defining the environment, a one-page code for the algorithm, and a refactored modularized neural network architecture.

In addition, we demonstrated that PPO is a viable algorithm for solving vehicle routing problems, as seen from its better performance with the Attention Model. Moving forward, our work can be further improved with more recent advances in OR models and RL algorithms. For instance, we only borrow the Multi-Greedy search strategies from POMO, yet its multi-start-nodes training strategy and inference augmentation strategy may also be helpful. We hope our work will serve as a new avenue for future vehicle routing problem research.

\section{Conclusion}
In this paper, we presented RLOR, a flexible framework of deep reinforcement learning for operation research that leverages recent advances in reinforcement learning algorithms. RLOR provides a comprehensive end-to-end deep learning model on a modern RL platform while improving the training efficiency of the original implementation of the Attention Model. We evaluated the performance of RLOR on the travelling salesman problem and capacitated vehicle routing problem. We demonstrated the improvement of RLOR over the original implementation of the Attention Model, showing the potential of the RLOR framework.

\bibliography{iclr2022_conference}
\bibliographystyle{iclr2022_conference}
\newpage
\appendix
\section{Algorithmic structures}\label{appendix-algo}
In this section, we will look at the algorithm structures of three works: PN+RL \citep{Bello2019NeuralLearning}, Attention Model \citep{Kool2019AttentionProblems}, and POMO \citep{Kwon2020POMO:Learning}. Without specification, the \texttt{train} function describes the algorithm in a \textbf{single} epoch.

\subsection{Basic algorithm} \label{sec:basic_algo}
\begin{minted}{python}
# REINFORCE Algorithm (basic)
def train(policy network @$p_\theta$@, training set @$S$@, batch size @$B$@):
    for @$i$@ in @$1...B$@:
         @$s_i$@ = RandomInstance(@$S$@)
         @$\pi_i$@ = SampleRollout(@$s_i,p_\theta$@)
    @$\vb$@ = UpdateBaseline(@$\vs,\bm{\pi}$@)
    @$\nabla\mathcal{L}$@ = @$\frac{1}{B}\sum_{i=1}^B (L(\pi_i|s_i)-b_i)\nabla_\theta \log(p_\theta(\pi_i|s_i))$@
    @$\theta$@ = GradientDescent(@$\theta,\nabla\mathcal{L}$@)
\end{minted}

Note that, the same REINFORCE algorithm can be used to fine-tune on a testing sample. It is denoted as Active Search in PN+RL from  \citet{Bello2019NeuralLearning}.

\begin{minted}{python}
# Fine tuning with Active Search
def finetune(policy network @$p_\theta$@, testing sample @$s$@, batch size @$B$@):
    for @$i$@ in @$1...B$@:
         @$\pi_i$@ = SampleRollout(@$s,p_\theta$@) # roll for a single instance
    @$\vb$@ = UpdateBaseline(@$s,\bm{\pi}$@)
    @$\nabla\mathcal{L}$@ = @$\frac{1}{B}\sum_{i=1}^B (L(\pi_i|s)-b_i)\nabla_\theta \log(p_\theta(\pi_i|s))$@
    @$\theta$@ = GradientDescent(@$\theta,\nabla\mathcal{L}$@)
\end{minted}

\subsection{Modifications}
We describe the changes different works have made.

\begin{minted}[highlightlines={3,8,11,12},highlightcolor=diffplus]{python}
 # Attention model : new baseline
 def train(policy network @$p_\theta$@, training set @$S$@, batch size @$B$@, 
+          significance @$\alpha$@):
     for @$i$@ in @$1...B$@:
         @$s_i$@ = RandomInstance(@$S$@)
         @$\pi_i$@ = SampleRollout(@$s_i,p_\theta$@)
@\Colorbox{diffminus}{-    $\vb$ = UpdateBaseline($\vs,\bm{\pi}$)}@
+    @$\vb$@ = UpdateBaseline(@$\vs,\bm{\pi},p_{\theta^\mathrm{BL}}$@)
     @$\nabla\mathcal{L}$@ = @$\frac{1}{B}\sum_{i=1}^B (L(\pi_i|s_i)-b_i)\nabla_\theta \log(p_\theta(\pi_i|s_i))$@
     @$\theta$@ = GradientDescent(@$\theta,\nabla\mathcal{L}$@)
+    if OneSidedPairedTest(@$p_\theta,p_{\theta^\mathrm{BL}}$@) < @$\alpha$@: # $p_\theta$ is better than $p_{\theta^\mathrm{BL}}$
+        @$\theta^\mathrm{BL}$@ = @$\theta$@
\end{minted}

\begin{minted}[highlightlines={3,7,8,11},highlightcolor=diffplus]{python}
 # POMO : new policy sampler
 def train(policy network @$p_\theta$@, training set @$S$@, batch size @$B$@, 
+          number of start nodes @$N$@):
    for @$i$@ in @$1...B$@:
         @$s_i$@ = RandomInstance(@$S$@)
@\Colorbox{diffminus}{-        $\pi_i$ = SampleRollout($s_i,p_\theta$)}@
+        @$\alpha_{i}^1,...,\alpha_{i}^N$@ = SelectStartNodes(@$s_i$@)
+        @$\pi_{i}^1,...,\pi_{i}^N$@ = SampleRollout(@$s_i,p_\theta,\{\alpha_{i,j}\}$@)
     @$\vb$@ = UpdateBaseline(@$\vs,\bm{\pi}$@)
@\Colorbox{diffminus}{-    $\nabla\mathcal{L}$ = $\frac{1}{B}\sum_{i=1}^B (L(\pi_i|s_i)-b_i)\nabla_\theta \log(p_\theta(\pi_i|s_i))$}@
+    @$\nabla\mathcal{L}$@ = @$\frac{1}{NB}\sum_{i=1}^B\sum_{j=1}^N (L(\pi_i^j|s_i)-b_i)\nabla_\theta \log(p_\theta(\pi_i^j|s_i))$@
     @$\theta$@ = GradientDescent(@$\theta,\nabla\mathcal{L}$@)
\end{minted}

\subsection{Baselines}
For different works, they use different baselines.
\begin{minted}{python}
# PN+RL
def UpdateBaseline(@$\vs,\bm{\pi}$@):
    for @$i$@ in @$1...B$@:
        @$b_i$@ = criticNetwork(@$s_i$@)
    return @$b$@
    
# Attention model
def UpdateBaseline(@$\vs,\bm{\pi},p_{\theta^\mathrm{BL}}$@):
    for @$i$@ in @$1...B$@:
        @$\pi_i^{\mathrm{BL}}$@ = GreedyRollout(@$s_i,p_{\theta^\mathrm{BL}}$@) # $p_{\theta^\mathrm{BL}}$ is the best model so far
        @$b_i$@ = @$L(\pi_i^{\mathrm{BL}}|s_i)$@
    return @$b$@
    
# POMO
def UpdateBaseline(@$\vs,\bm{\pi}$@):
    for @$i$@ in @$1...B$@:
        @$b_i$@ = @$\frac{1}{N}\sum_{j=1}^N(L(\pi_i^j|s_i))$@ # $N$ rollouts for each $\pi_i$
    return @$b$@
\end{minted}

\subsection{Decoding strategies}
In this subsection, we denote $\pi_t$ to be the action chosen at step $t$.

\begin{minted}{python}
# Greedy strategy
def GreedyRollout(instance @$s$@, policy network @$p_\theta$@):
    @$\pi$@ = []
    while not Done(@$s,\pi$@):
        @$\pi_{next}$@ = @$\mathrm{argmax}_i\,{ p_\theta(i|s,\pi)}$@
        @$\pi$@.append(@$\pi_{next}$@)
    return @$\pi$@
    
def inference(instance @$s$@, policy network @$p_\theta$@):
    @$\pi$@ = GreedyRollout(@$s$@,@$p_\theta$@)
    return @$\pi$@
\end{minted}
\begin{minted}{python}
# Sampling strategy
def SampleRollout(instance @$s$@, policy network @$p_\theta$@):
    @$\pi$@ = []
    while not Done(@$s,\pi$@):
        @$\pi_{next}$@ = Sample(@${ p_\theta( . |s,\pi)}$@)
        @$\pi$@.append(@$\pi_{next}$@)
    return @$\pi$@
    
def inference(instance @$s$@, policy network @$p_\theta$@, sampling size @$N$@):
    for @$i$@ in @$1...N$@:
        @$\pi_i$@ = SampleRollout(@$s$@,@$p_\theta$@)
    @$\pi_{bext}$@ = @$\mathrm{argmin}_{\pi_i}\,{ L(\pi_i|s)}$@
    return @$\pi_{best}$@
\end{minted}
\begin{minted}{python}
# Beam Search
def inference(instance @$s$@, policy network @$p_\theta$@, beam size @$w$@):
    beam = [([],0)] # [seq, score = $\sum_{\pi_t} \log p_\theta(\pi_t|s,\pi_{1:t-1})$]
    finalBeam = []
    while beam is not empty:
        expansion = []
        for (@$\pi$@, score) in beam:
            if Done(@$s$@,@$\pi$@):
                finalBeam.append(@$\pi$@)
                continue
            for j in @$p_\theta(.|s,\pi)_{>0}$@: # only feasible actions
                expansion.append((@$\pi$@+[j], score + @$\log p_\theta(\pi_t=j|s,\pi)$@))
        beam = TopK(expansion)
    @$\pi_{bext}$@ = @$\mathrm{argmin}_{\pi\in \texttt{finalBeam}}\,{ L(\pi|s)}$@
    return @$\pi_{bext}$@
\end{minted}

As mentioned in Section \ref{sec:basic_algo}, the model can be fine-tuned during inference with Active Search.
\begin{minted}{python}
# Active Search
def finetune(policy network @$p_\theta$@, testing sample @$s$@, batch size @$B$@):
    for @$i$@ in @$1...B$@:
         @$\pi_i$@ = SampleRollout(@$s,p_\theta$@)
    @$\vb$@ = UpdateBaseline(@$s,\bm{\pi}$@)
    @$\nabla\mathcal{L}$@ = @$\frac{1}{B}\sum_{i=1}^B (L(\pi_i|s)-b_i)\nabla_\theta \log(p_\theta(\pi_i|s))$@
    @$\theta$@ = GradientDescent(@$\theta,\nabla\mathcal{L}$@)

def inference(policy network @$p_\theta$@, testing sample @$s$@, batch size @$B$@,
              epochs @$T$@):
    for @$t$@ in @$1...T$@:
        finetune(@$p_\theta$@,@$s$@,@$B$@)
    @$\pi_{bext}$@ = best policy encounted during fine-tuning
    return @$\pi_{bext}$@ 
    
\end{minted}

We can integrate the end2end learning model with classical exact algorithms and use our policy network output as a heuristic. This converts the end2end learning model to an exact model (we can run exhaustive search under the classical exact algorithm framework). 
\begin{minted}{python}
def DepthFirstBranchAndBound(instance @$s$@, policy network @$p_\theta$@ ):
    bound = @$\infty$@
    solution = None
    def dfs(@$\pi$@):
        # bound stage: prune
        # Prune if the partial solution is worse than the incumbent 
        if @$L(\pi|s)$@ >= bound: 
            return
        # bound stage: update bound
        # Update the bound if a better solution is found
        if Done(@$s,\pi$@) and @$L(\pi|s)$@ < bound:
            bound = @$L(\pi|s)$@
            solution = @$\pi$@
            return
 
        # branch stage: RL prediction as the heuristic
        values = @$\mathrm{argsort}_j p_\theta(\pi_t=j|s,\pi)$@
        for action in values: 
            @$\pi_{new}$@ = copy(@$\pi$@).append(action)
            dfs(@$\pi_{new}$@)
    @$\pi_{0}$@ = []
    dfs(@$\pi_{0}$@)
    return solution
\end{minted}

\newpage
\section{Models}\label{appendix-models}
In this section, we will look at the neural network architecture of two works: Pointer Network \citep{Vinyals2015PointerNetworks}, and Attention Model \citep{Kool2019AttentionProblems}.

\subsection{Pointer Network}
\begin{minted}{python}
Embedding(@$s$@) 
    -> @$\vx=W_0s$@

Encoder(@$\vx$@) # 1 layer lstm
    -> @$\vo,(h,c)$@ = LSTM(@$\vx$@)  # output, (hidden state, cell state)
    
Attention(@$q,\vk$@) 
    -> (@$W_1k_j$@, @$u_j=v^T \mathrm{tanh}(W_1k_j+W_2q)$@) @$\text{for each key }k_j$@
    
Attention@$_{clip}$@(@$q,\vk$@) # clip the range of logits to [-C,C]
    -> (@$W_1k_j$@, @$C\mathrm{tanh}(u_j)$@ )
    
Decoder(@$\vx,(h_0 ,c_0), \vo$@):
    for @$t$@ in 1...T:
        # Predict Step
        @$h_t, c_t=\mathrm{LSTMCell}(\vx_{\pi_{t-1}},(h_{t-1},c_{t-1}))$@
        @$\ve_t, \vu_t=\mathrm{Attention}(h_t,  \vo)$@
        @$\va_t = \mathrm{Softmax}(\mathrm{Masked}_t(\vu_t))$@
        @$h_t^{\prime}=\sum_i a_{t,i}e_{t,i}$@ # glimpse
        @$\ve_t^{\prime}, \vu_t^{\prime}=\mathrm{Attention_{clip}}(h_t^{\prime},  \vo)$@
        @$\vp_t = \mathrm{Softmax}(\mathrm{Masked}_t(\vu_t^{\prime}))$@
        # Decode Step
        @$\pi_{t} = \mathrm{DecodingStartegy}(\vp_t)$@
        # Update Step
        @$\mathrm{Masked}_{t+1} = \mathrm{Masked}_t.update(\pi_t)$@
    return @$\{\log(\vp_t)\},\pi$@
    
PointerNetwork(@$s$@):
    @$\vx$@ = Embedding(@$s$@)
    @$\vo,(h,c)$@ = Encoder(@$\vx$@) 
    @$\{\log(\vp_t)\},\pi$@ = Decoder(@$\vx,(h_{-1} ,c_{-1}), \vo$@)
    @$\log(p_\theta(\pi|s))$@ = @$\sum_t\log(\vp_{t,\pi_t})$@
    return @$\log(p_\theta(\pi|s)),\pi$@
\end{minted}
\newpage

\subsection{Attention Model} \label{appendix-attention-model}
\begin{minted}{python}
#####################################
###  Helper classes for attention ### 
#####################################

AttentionScore(@$q,\vk, \mathrm{mask}$@) 
    -> @$u_j=\begin{cases}\frac{q^Tk_j}{\sqrt{\smash{d_q}}}\text{ for node }j \notin \mathrm{mask}\\-\infty  \text{ otherwise. }\end{cases} $@
    
AttentionScore@$_{clip}$@(@$q,\vk, \mathrm{mask}$@)  # clip the range of logits to [-C,C]
    -> @$u_j=\begin{cases}C\mathrm{tanh}(\frac{q^Tk_j}{\sqrt{\smash{d_q}}})\text{ for node }j \notin \mathrm{mask}\\-\infty  \text{ otherwise. }\end{cases} $@

MultiHeadAttention(@$q,\vk,\vv, \mathrm{mask}$@):
    @$\pmb{a}^{(j)} = \mathrm{Softmax}(\mathrm{AttentionScore}(q^{(j)},\pmb{k}^{(j)}, \mathrm{mask})) \text{ for each head }j$@
    @$h^{(j)} = \sum\nolimits_i \pmb{a}^{(j)}_i\pmb{v}_i$@  #reweigh value with attention weight in each head
    @$q^\prime = W^O \left[h^{(1)},...,h^{(J)}\right]$@  # flatten to single head
    return @$q^\prime$@
    
MultiHeadAttentionProj(@$q_0,\vh, \mathrm{mask}$@): # with projection
    @$q, \pmb{k}, \pmb{v} = W^Qq_0, W^K\vh, W^V\vh$@
    @$q^\prime = \mathrm{MultiHeadAttention}(q,\vk,\vv, \mathrm{mask})$@
    return @$q^\prime$@
\end{minted}
\newpage
\begin{minted}{python} 
#####################################
### Main code for attention model ### 
#####################################

Embedding(@$s$@) 
    -> @$\vx=\left[W_0[s_{loc},s_{fea}], W_{depot}s_{depot}\right]$@
    
MHALayer(@$\pmb{x}$@): # self attention, then MLP
        @$\pmb{x}_0=\pmb{x}+\mathrm{MultiHeadAttentionProj}(\pmb{x})$@
        @$\pmb{x}_1=\mathrm{BatchNorm}(\pmb{x}_0)$@
        @$\pmb{x}_2=\pmb{x}_1+\mathrm{MLP_{\text{2 layers}}}(\pmb{x}_1)$@
        @$\pmb{h} =\mathrm{BatchNorm}(\pmb{x}_2)$@
        return @$\pmb{h}$@
        
Encoder(@$\vx$@) # 2 MHA layers
    @$\pmb{h} = \pmb{x}$@ 
    for i in n_layers:
        @$\pmb{h} = \mathrm{MHALayer_i}(\pmb{h})$@ 
    return @$\pmb{h}$@ 

Context(@$s,\vh$@):
    -> @$(\vh_{\pi_{t-1}}, \vh_{s}) = \text{(previous node, global state (e.g. remaining capacity) }$@

Decoder(@$\vx,\vh$@):
    for @$t$@ in 1...T:
        # Predict Step
        @$\bar{\pmb{h}} = \frac{1}{N}\sum\nolimits_i \pmb{h}_i$@  # graph embedding
        @$\pmb{h}_{(c)} = [\bar{\pmb{h}}, \text{Context}(s, \pmb{h})]$@
        @$q_{gl} = \mathrm{MultiHeadAttention}(W^Q \pmb{h}_{(c)},W^K \pmb{h},W^V \pmb{h},\mathrm{mask}_t)$@
        @$\pmb{p}_t = \mathrm{Softmax}(\mathrm{AttentionScore}_{\text{clip}}(q_{gl},W^{K^\prime} \pmb{h}, \mathrm{mask}_t))$@
        # Decode Step
        @$\pi_{t} = \mathrm{DecodingStartegy}(\vp_t)$@
        # Update Step
        @$\mathrm{mask}_{t+1} = \mathrm{mask}_t.update(\pi_t)$@
    return @$\{\log(\vp_t)\},\pi$@
    
AttentionModel(@$s$@):
    @$\vx$@ = Embedding(@$s$@)
    @$\vh$@ = Encoder(@$\vx$@) 
    @$\{\log(\vp_t)\},\pi$@ = Decoder(@$s,\vh$@)
    @$\log(p_\theta(\pi|s))$@ = @$\sum_t\log(\vp_{t,\pi_t})$@
    return @$\log(p_\theta(\pi|s)),\pi$@
\end{minted}
Note that there are some tricks to speed up the training and inference for the decoder, such as precomputing and caching the projection of the graph embedding $\bar{\pmb{h}}$, as well as the keys and values $W^K \pmb{h}, W^V \pmb{h}, W^{K^\prime} \pmb{h}$. 

For dynamic node features, their projections are recomputed in every decoding steps with $\mathrm{DynamicEmbedding}$. The projections are added to the keys and values to form the new keys and values in the decoder.
\newpage
\section{More on Vectorized Environment}\label{appendix-vectorization}
A vectorized environment runs multiple environments, either sequentially/ in parallel/ in other way. It receives a batch of actions, runs environment steps, and returns a batch of observations at that step. The way a RL library manipulates multiple environments (including creation, \textit{step}, \textit{reset}) is denoted as the vector API.

Among the RL platforms, CleanRL has the most efficient vector API --- the one from OpenAI Gym, instead of writing its own. RLlib is also planning to depreciate its own vector API and migrate it to Gym's vector API as indicated in the \href{https://github.com/Farama-Foundation/Gymnasium/issues/32#issuecomment-1357945491}{github repo}.

The Gym's vector API is a gym environment. It can be treated as a wrapper to interact with multiple environments. Take \mintinline{python}{gym.vector.SyncVectorEnv} as an example. When it is instantiated, it sequentially instantiates each of the sub-environments. In each \mintinline{python}{step}, it sequentially feeds the actions to the sub-environments, receives the results from them and stores the \mintinline{python}{reward}, \mintinline{python}{done} directly to the pre-allocated numpy arrays while the \mintinline{python}{observations} are concatenated according to their type (dict, list, tuple, float, int, ... ). To handle the concatenation/ array pre-allocation of different types, Gym's vector API makes heavy use of single-dispatch function instead of nested \mintinline{python}{isinstance} as in DI-engine. The single-dispatch function approach is cleaner and is easier to extend.

We compared \mintinline{python}{SyncVectorEnv} (run sequentially) and \mintinline{python}{AsyncVectorEnv} (run in parallel) on our routing problem environments. We found that the sequential one is faster. It may be attributed to the communication overhead across the pipes for the parallel settings.

\paragraph{Explicitly vectorized environment.} Instead of using Gym's vector API, we can implement the vectorized environment explicitly. For example, the environments in the Attention Model paper \citep{Kool2019AttentionProblems} are written with PyTorch tensors and PyTorch's vectorized operations. The Gym environment for a routing problem written with numpy can also be re-implemented with numpy's vectorized operations. The explicit implementation of vectorized environment provides the same interface as the vector API, yet allows more room of optimization in computation efficiency. 

\subsection{Chunking Vectorized Environment}
Thanks to the flexibility of the Gym's vector API, it is possible to implement chunking: vectorization of multiple vectorized envinronment. There are a couple of advantages of chunking: efficiency and flexibility. The explicitly vectorized environments are faster than running individual environments sequentially/ in parallel with the vector API. However, it requires careful design for fine-grained control in these explicitly vectorized environments. With chunking, we do not need to dive into the implementation of explicitly vectorized environments. Instead, we can apply different wrappers to these vectorized environments. 

An example usage of chunking is that, we want to create $M$ TSP problems and explore $N$ trajectories in each of the TSP problems. First, we can create a vectorized environment (denote as $\mathrm{TrajEnv}_i$) of $N$ sub-environments. We can apply a wrapper to ensure the same nodes coordinates in these sub-environments. Then, we create a higher level of vectorized environment (denote as $\mathrm{ProblemEnv}$) of these $M$ sub-environments $\{\mathrm{TrajEnv}_i, i = 1...M\}$. The bi-level vectorized environment can now well describe the setting we required. 

With Gym's vector API, the \mintinline{python}{observations} from the bi-level vectorized environment has its dimensions extended accordingly. It provides clearer semantics for performing operations on observations from different level of vectorized environments. For example, the static features from the same $\mathrm{TrajEnv}_i$ need only to be computed once and can be shared among different trajectories in the sub-environments. It is straight-forward to do so with Gym's vector API, otherwise it would require careful reshaping and indexing.

\end{document}